\documentclass[12pt]{article}
\usepackage[latin1]{inputenc}
\usepackage{amssymb,amsmath,amsfonts}
\newtheorem{theorem}{Theorem}
\newtheorem{proposition}[theorem]{Proposition}
\newtheorem{lemma}[theorem]{Lemma}

\newtheorem{remark}[theorem]{Remark}

\newtheorem{example}[theorem]{Example}

\newcommand{\F}{\mathcal{F}}
\newcommand{\R}{\mathbb{R}}
\newcommand{\Q}{\mathbb{Q}}
\newcommand{\Sf}{\mathbb{S}}
\newcommand{\Hy}{\mathbb{H}}
\renewcommand{\dim}{\mbox{dim\,}}
\newcommand{\spa}{\mbox{span\,}}
\newcommand{\hess}{\mbox{Hess\,}}
\newcommand{\Ric}{\mbox{Ric}}

\newcommand{\grad}{\mbox{grad\,}}
\newcommand{\po}{{\hspace*{-1ex}}{\bf .  }}

\newcommand{\Lo}{\mathbb{L}}
\newcommand{\V}{\mathbb{V}}
\newcommand{\E}{\mathbb{E}}
\def\<{\langle}
\def\>{\rangle}
\def\n{\nabla}
\def\d{\partial}
\def\a{\alpha}
\def\b{\beta}
\def\be{\begin{equation} }
\def\ee{\end{equation} }

\def\nap{\nabla^\perp}
\def\proof{\noindent{\it Proof:  }}
\def\qed{\ifhmode\unskip\nobreak\fi\ifmmode\ifinner
\else\hskip5 pt \fi\fi\hbox{\hskip5 pt \vrule width4 pt
height6 pt  depth1.5 pt \hskip 1pt }}

\begin{document}

\title{Conformally flat submanifolds with flat\\ normal bundle}
\author{M.\ Dajczer, C.-R.\ Onti and Th.\ Vlachos}
\date{}
\maketitle
\vspace*{-4ex}
\begin{center}
In memory of Manfredo do Carmo 
\end{center}

\begin{abstract} 
We prove that any conformally flat submanifold with  flat normal bundle 
in a conformally flat Riemannian manifold is locally holonomic, that is, 
admits a principal coordinate system. As one of the consequences of this 
fact, it is shown that the Ribaucour transformation can be used to construct 
an associated large family of immersions with induced conformally flat metrics 
holonomic with respect to the same coordinate system.
\end{abstract}

A main task in conformal geometry is the study of submanifolds of conformally 
flat Riemannian manifolds with induced conformally flat metrics. A Riemannian 
manifold $M^n$ is said to be \emph{conformally flat} if each point lies in an 
open neighborhood conformal to an open subset of Euclidean space $\R^n$. 
This is always the case for manifolds endowed with metrics of constant sectional 
curvature.

Even if they belong to the realm of conformal geometry, for reason of 
simplicity most of the results in this paper are stated for submanifolds 
of Euclidean space. Nevertheless, they hold true when the ambient space is 
just a conformally flat manifold. 

E. Cartan \cite{ca0} proved  that a hypersurface $f\colon M^n\to\R^{n+1}$, 
$n\geq 4$,  is conformally flat if and only if at each point there is a 
principal curvature of multiplicity at least $n-1$.  If $M^n$ is free of flat 
points, then $f(M)$  is locally foliated by $(n-1)$-dimensional umbilical 
submanifolds of $\R^{n+1}$, or equivalently, we have that  $f(M)$ is enveloped 
by a one-parameter family of umbilical hypersurfaces of the ambient space.

Moore \cite{mo0} extended Cartan's result to submanifolds of higher 
codimension. He showed that an isometric immersion of a  conformally flat 
manifold $f\colon M^n\to\R^{n+p}$ of dimension  $n\geq 4$ and codimension 
$p\leq n-3$ has a principal normal vector of multiplicity at least $n-p\geq 3$ 
at each point. Recall that a normal vector $\eta\in N_fM(x)$  is called a 
\emph{principal normal} of $f$ at $x\in M^n$ with \emph{multiplicity} $s$ 
if the tangent subspace defined as
$$
E_\eta(x)=\big\{X\in T_xM:\alpha_f(X,Y)=\<X,Y\>\eta\;\;
\text{for all}\;\;Y\in T_xM\big\}
$$
in terms of the second fundamental form
$\alpha_f\colon TM\times TM\to N_f M$ of the immersion,
satisfies  $\dim E_\eta(x)=s>0$. Clearly, principal normals are a natural
generalization to submanifolds of higher codimension of principal curvatures
of hypersurfaces.

A smooth normal vector field $\eta\in N_fM$ to an isometric immersion
$f\colon M^n\to\R^N$ is called a \emph{principal normal vector field} 
with \emph{multiplicity} $s$ if $\dim E_\eta(x)=s>0$ is constant, in 
which case the distribution $x\in M^n\mapsto E_\eta(x)$ is smooth. 
A principal normal vector field $\eta$ is called \emph{Dupin} if it is 
parallel along $E_\eta$ in the normal connection, which is always the 
case if $s\geq 2$.  The principal normal being Dupin implies that $f$ 
maps each leaf of the spherical distribution $x\mapsto E_\eta(x)$ into 
an umbilical submanifold of $\R^{n+p}$.

Let $f\colon M^n\to\R^N$ be an isometric immersion with flat normal bundle, 
that is, at any point the curvature tensor of the metric induced from the 
ambient space on the normal bundle of the submanifold vanishes. Submanifolds 
with flat normal bundle have  captured the attention because they 
``behave like hypersurfaces". 
For instance, from \cite{re} we have that at any point $x\in M^n$ there 
is a  unique set of pairwise distinct principal normal vectors 
$\eta_i\in N_fM(x), 1\leq i\leq s(x)$, and an associate orthogonal 
splitting of the tangent space as
$$
T_xM=E_{\eta_1}(x)\oplus\cdots\oplus E_{\eta_s}(x).
$$
Related to Moore's result in \cite{mo0} it was proved in \cite{dt2}
that if at some point a conformally flat submanifold 
$f\colon M^n\to\R^{2n-2}$, $n\geq 4$, has no principal normal of
multiplicity larger than one, then the normal bundle at that point 
has to be flat.

An isometric immersion $f\colon M^n\to\R^N$  with flat normal bundle 
is called  \emph{holonomic} if $M^n$ carries a global orthogonal principal 
coordinate system. That the coordinates are principal means that the 
corresponding coordinate vector fields diagonalize the second fundamental 
form of the immersion at any point. It is a classical fact that the 
Gauss-Codazzi equations for a holonomic submanifold can be nicely 
written as a completely integrable system of first order PDE's.

Cartan \cite{ca} proved that if $f\colon M^n_{\tilde c}\to\Q_c^{n+p}$ is 
an isometric immersion of a manifold with constant sectional curvature 
$\tilde{c}$ into a space form of sectional curvature $c$, then the 
submanifold is locally holonomic if $\tilde{c}<c$ and the codimension 
is $p=n-1$, which in this case is the least possible. 
If $\tilde{c}>c$  the same conclusion for the same codimension was 
obtained by Moore \cite{mo}  under the additional assumption that 
the submanifold is free of weak-umbilic points. To prove these results, 
one first has to argue that the submanifold must have flat normal bundle, 
and then that the image of its second fundamental form spans at any point
the full normal space of the immersion. An elementary argument then yields
local holonomicity; for instance see Proposition $1$ in \cite{dt}. 

Holonomic isometric immersions $f\colon M^n_{\tilde c}\to\Q_c^N$ are of 
particular interest because the associated Gauss-Codazzi system of equations 
is, in this case, a natural generalization of the sinh-Gordon, sine-Gordon, 
Laplace or wave equation, according to the values of the sectional 
curvatures. These equations are classically known to be in correspondence 
to constant curvature surfaces; cf.\ \cite{dt0} and \cite{dt}.

Related to the above, we proved in \cite{dov} that any proper isometric 
immersion of an Einstein manifold $f\colon M^n\to\R^N$ with flat normal 
bundle  is locally holonomic. Throughout this paper that a submanifold 
with flat normal bundle is proper means that it has a constant number 
of principal normals. 
Note that the only Riemannian manifolds that are simultaneously 
Einstein and conformally flat are the ones of constant sectional curvature.
\vspace{1ex}

The following is the main result of this paper.

\begin{theorem}\po\label{one} Let $f\colon M^n\to\R^N,n\geq 4$, be a proper
isometric immersion with flat normal bundle of a conformally 
flat manifold. Then $f$ is locally holonomic with at most one 
principal normal vector field of multiplicity larger than one.  
\end{theorem}

In Example \ref{ex} given below a large family of nonflat conformally flat 
submanifolds of codimension two with flat normal bundle is constructed. 
They possess three principal normal vector fields and the holonomic coordinates 
are provided by the construction.
\vspace{1ex}

In view of the above result, it is quite natural to consider that the
conformally flat proper submanifolds with flat normal bundle 
belong to one of two classes according to whether they carry a principal 
normal vector field of multiplicity larger than one or not. Notice that in low 
codimension, the former situation is always the case due to the 
aforementioned result of Moore. 
\vspace{1ex}

An isometric immersion $f\colon M^n\to\R^{n+p}$ is said to be 
\emph{quasiumbilical} if at any point of $M^n$ there exists an orthonormal 
normal base $\xi_1,\ldots,\xi_p$ such that each 
shape operator $A_{\xi_j}$, $1\leq j\leq p$, has an eigenvalue of 
multiplicity at least $n-1$. The property of a submanifold being 
quasiumbilical is conformally invariant.
\vspace{1ex}

In view of Cartan's result that isometric immersions with codimension 
one between conformally flat manifolds  $f\colon M^n\to\tilde M^{n+1}$, 
$n\geq 4$, are quasiumbilical, it is clear that composing hypersurfaces 
of this type yields quasiumbilical submanifolds with higher codimension.

If $f\colon M^n\to\R^{n+p}$, $n\geq 4$, is quasiumbilical, it is easy to 
see that the Weyl tensor of $M^n$ vanishes, hence $M^n$ is conformally flat. 
On the other hand, Chen and Verstraelen \cite{bycv} proved that if $M^n$, 
$n\geq 4$, is conformally flat and $f\colon M^n\to\R^{n+p}$ has flat normal 
bundle with  codimension $p\leq n-3$, then the submanifold is quasiumbilical. 
Moore and Morvan \cite{momo} reached the same conclusion 
for codimension $p\leq 4$ without the assumption of flatness of the normal 
bundle. 

To conclude that a conformally flat submanifold with flat normal bundle 
is quasiumbilical, the presence for a  principal normal of multiplicity at 
least two suffices regardless of the codimension. 

\begin{theorem}\po\label{quasi}
Let $f\colon M^n\to\R^{n+p}$, $n\geq 4,$ be a proper isometric immersion 
with flat normal bundle of a conformally flat manifold. 
If $f$ carries a principal normal vector field of multiplicity 
$m\geq 2$ then $p\geq n-m$ and $f$ is quasiumbilical.
\end{theorem}

If a principal normal is trivial, say $\eta_i=0$, then  
$$
E_{\eta_i}(x)=\Delta(x)=\{X\in T_xM:\alpha_f(X,Y)=0
\;\;\mbox{for all}\;\;Y\in T_xM\}
$$
is called the \emph{relative nullity} subspace of $f$ at $x\in M^n$ and 
$\nu(x)=\dim\Delta(x)$ the \emph{index of relative nullity} of $f$ 
at $x\in M^n$.
\vspace{1ex}

That a conformally flat submanifold with flat normal bundle has index $\nu\geq 1$
turns out to be quite restrictive. It is convenient to state the following 
result for ambient space forms $\Q_c^N$ and leave the definitions of 
generalized cone and cylinder in these spaces for later. 

\begin{theorem}\label{gen}\po
Let $f\colon M^n\to\Q_c^N, n\geq 4$, be a proper isometric immersion with 
flat normal bundle of a conformally flat manifold.
If $f$  has index of relative nullity $\nu\geq 1$,
then one of the following holds: 
\begin{itemize}
\item[(i)] $M^n$ has constant sectional curvature $c$ and $f$ is locally 
a $\nu$-generalized cylinder over a holonomic submanifold 
$g\colon L^{n-\nu}\to\Q_c^N$. 
\item[(ii)] $M^n$ has sectional curvature different from $c$ and 
$f$ is locally  a $1$-generalized cone over a holonomic submanifold 
$g\colon L^{n-1}\to \Q_{\tilde c}^{N-1}\subset\Q_c^N$, $\tilde{c}\geq c$, 
with constant sectional curvature different from $\tilde c$.
\end{itemize}
\end{theorem}

Holonomic submanifolds are the natural object of application of the 
Ribaucour transformation introduced in \cite{dt}. This fact is 
instrumental to obtain the following result.

\begin{theorem}\po\label{two} Let $f\colon M^n\to\R^N$, $n\geq 4$, 
be a proper isometric immersion with flat normal bundle of a 
conformally flat  manifold. Then locally there exists an 
$N$-parameter family of immersions $\tilde{f}\colon M^n\to\R^N$ with 
induced conformally flat metrics that are holonomic with respect to 
the same coordinate system as $f$.
\end{theorem}

We observe that some of the results in this paper have been obtained 
by Donaldson and Terng \cite{nt} under strong additional assumptions.

\section{Preliminaries}

In this section we show that the statements in this paper
are conformally invariant.\vspace{2ex}

Let $f\colon M^n\to\R^N$ be an isometric immersion with flat 
normal bundle and let $\eta_i\in N_fM(x)$, $1\leq i\leq s(x)$, 
be the set of pairwise distinct principal normals at $x\in M^n$. 
Then, the second fundamental form $\alpha=\alpha_f$ of $f$ 
acquires the form
$$
\alpha(X,Y)(x)=\sum_{i=1}^s\<X^i,Y^i\>\eta_i
$$
where $X\mapsto X^i$ denotes the orthogonal 
projection from $T_xM$ onto $E_i(x)=E_{\eta_i}(x)$.
Equivalently, in terms of the shape operators of $f$ we have 
\be\label{equiv}
A_\xi X=\sum_{i=1}^s\<\xi,\eta_i\>X^i
\ee
for any $\xi\in N_fM$.
\vspace{1ex}

A submanifold $f\colon M^n\to\R^N$ with flat normal bundle is 
called \emph{proper} if $s(x)=k$ is constant on $M^n$. In this situation, 
we have from \cite{re} that the principal normal vector fields 
$x\in M^n\mapsto \eta_i(x)$, $1\leq i\leq k$, are smooth. Moreover, the  
distributions $x\in M^n\mapsto E_i(x)$, $1\leq i\leq k$, have constant 
dimension and are also smooth. 
\vspace{1ex}

 Let $\tilde M^N$ be endowed with conformal metrics $g_1$ and $g_2$, 
that is, $g_2=\lambda^2g_1$ where $\lambda\in C^\infty(\tilde M)$ is 
positive. Given an immersion $f\colon M^n\to\tilde M^N$,  we thus have 
the two isometric immersions with the induced metrics
$$
f_j=f\colon (M^n, f^*g_j)\to (\tilde M^N, g_j),\;\;1\leq j\leq 2.
$$
At any $x\in M^n$  the  second fundamental forms of $f_1$ 
and $f_2$ are related by
$$
\alpha_{f_2}(X,Y)=\alpha_{f_1}(X,Y)
-\frac{1}{\lambda}g_1(X,Y)(\grad_1\lambda)^\perp
$$
and the normal curvature tensors by
$$
R^\perp_2(X,Y)\xi=R^\perp_1(X,Y)\xi
$$
for any $X,Y\in T_xM$ and $\xi\in N_fM(x)$.
In particular, if $\eta$ is a principal normal vector of $f_1$ 
at $x\in M^n$ then 
\be\label{corres}
\eta-\frac{1}{\lambda}(\grad_1\lambda)^\perp
\ee
is a principal normal vector of $f_2$ at $x\in M^n$. 
\vspace{1ex}

We thus have the following fact.

\begin{proposition}\po\label{metrics}
Let $f\colon M^n\to M^N$ be a proper isometric immersion with flat 
normal bundle and let $\tau\colon M^N\to\tilde{M}^N$ 
be a conformal diffeomorphism. Then the conformal immersion
$\tilde f=\tau\circ f\colon M^n\to\tilde{M}^N$ also has
flat normal bundle and is proper. 
\end{proposition}

\section{Proof of Theorem \ref{one}}

The proof of Theorem \ref{one} will follow from the two 
lemmas given in the sequel.

\begin{lemma}\po\label{unique} Let $f\colon M^n\to\R^N$, $n\geq 4$, 
be an isometric immersion with flat normal bundle of a conformally 
flat manifold. Then at any point of $M^n$ there exists at most one 
principal normal of  multiplicity at least two.
\end{lemma}

\proof It is well-known that the curvature tensor of $M^n$ has 
the form
$$
R(X,Y,Z,W)=L(X,W)\<Y,Z\>-L(X,Z)\<Y,W\>+L(Y,Z)\<X,W\>-L(Y,W)\<X,Z\>
$$
in terms of the Schouten tensor given by
$$
L(X,Y)=\frac{1}{n-2}\left(\Ric(X,Y)-\frac{\tau}{2(n-1)}\<X,Y\>\right)
$$
where $\tau$ denotes the scalar curvature. In particular, 
the sectional curvature is given by
\be\label{seccur}
K(X,Y)=L(X,X)+L(Y,Y)
\ee
where  $X,Y\in TM$ are orthonormal vectors. 

A straightforward computation of the Ricci tensor using the 
Gauss equation
\be\label{eqgauss}
R(X,Y,Z,W)=\<\alpha(X,W),\alpha(Y,Z)\>-\<\alpha(X,Z),\alpha(Y,W)\>
\ee 
yields
\be\label{ricci}
\Ric(X,Y)
=n\<\alpha(X,Y),H\>-\sum_{j=1}^n\<\alpha(X,X_j),\alpha(Y,X_j)\>
\ee
where $H$ is the mean curvature vector and $X_1,\dots,X_n$ an 
orthonormal tangent basis. 

We obtain from \eqref{seccur} and \eqref{eqgauss} that
\be\label{dec2}
L(X,X)+L(Y,Y)=\<\a(X,X),\a(Y,Y)\>-\|\a(X,Y)\|^2
\ee
for any pair $X,Y\in TM$ of orthonormal vectors. From \eqref{ricci} we have 
$$
\Ric(X,X)=n\<\eta_i,H\>-\|\eta_i\|^2
$$
for any unit vector $X\in E_i$.  Thus
\be\label{same}
(n-2)L(X,X)=n\<\eta_i,H\>
-\|\eta_i\|^2-\frac{\tau}{2(n-1)}
\ee
for any unit vector $X\in E_i$. Denoting 
$$
\hat\eta_i=\eta_i-H,\;\; 1\leq i\leq k,
$$
we obtain from \eqref{same} that
\be\label{newsame}
(n-2)L(X,X)=(n-1)\|H\|^2+(n-2)\<\hat\eta_i,H\>
-\|\hat\eta_i\|^2-\frac{\tau}{2(n-1)}
\ee
for any unit vector $X\in E_i$.

Assume that $\eta_1$ has multiplicity at least two. 
Then \eqref{same} yields 
$$
L(X,X)=L(Y,Y)
$$
for any unit vectors $X,Y\in E_1$. Hence \eqref{dec2} gives that
\be\label{twol}
2L(X,X)=\|\hat\eta_1\|^2+2\<\hat\eta_1,H\>+\|H\|^2
\ee
for any unit vector $X\in E_1$. It follows from  \eqref{newsame} and  
\eqref{twol} that
\be\label{hat}
\|\hat\eta_1\|^2=\|H\|^2-\frac{\tau}{n(n-1)}\cdot
\ee
We obtain from \eqref{twol}  and \eqref{hat} that
\be\label{same2}
L(X,X)=\|H\|^2+\<\hat\eta_1,H\>-\frac{\tau}{2n(n-1)}
\ee
for any unit vector $X\in E_1$. 

Given principal normals $\eta_i\neq\eta_j$, we have from 
\eqref{dec2} that
\be\label{what}
L(X,X)+L(Y,Y)=\<\hat\eta_i,\hat\eta_j\>+\<\hat\eta_i+\hat\eta_j,H\>+\|H\|^2
\ee
where $X\in E_i$ and $Y\in E_j$ are unit vectors.
Suppose that $\eta_i$ and $\eta_j$ have both multiplicity at least two. 
It follows from \eqref{hat} that
\be\label{same4}
\|\hat\eta_i\|^2=\|H\|^2-\frac{\tau}{n(n-1)}=\|\hat\eta_j\|^2.
\ee
On the other hand, we obtain from \eqref{same2} and \eqref{what} that
\be\label{same3}
\<\hat\eta_i,\hat\eta_j\>=\|H\|^2-\frac{\tau}{n(n-1)}\cdot
\ee
We conclude from \eqref{same4} and \eqref{same3} that 
$\eta_i=\eta_j$, and this is a contradiction.\qed

\begin{example}\po {\em A rather simple example in high codimension 
of a conformally flat submanifold with flat normal bundle carrying a 
principal normal of multiplicity at least two is as follows: Let $M^{2n}$ 
be the Riemannian product $\Sf^n_1\times U$ where $\Sf^n_1\subset\R^{n+1}$ 
is a round sphere and $U$  an open subset of the hyperbolic space 
$\Hy_{-1}^n$ isometrically immersed in $\R^{2n-1}$. Then $M^{2n}$ is
conformally flat and the product isometric immersion of $M^{2n}$ into 
$\R^{3n}$ has a principal normal of multiplicity $n$.}
\end{example}

\begin{lemma}\label{li}\po Let $f\colon M^n\to\R^N, n \geq 4,$ be an 
isometric immersion with flat normal bundle of a conformally 
flat manifold. If at some point of $M^n$ we have  $k\geq 3$, then 
the vectors $\eta_j-\eta_m$ and $\eta_j-\eta_\ell$ are linearly 
independent for $1\leq m\neq j\neq \ell\neq m\leq k$.
\end{lemma}

\proof We argue by contradiction. In the sequel suppose that
\be\label{hipo}
\eta_j-\eta_m=\mu(\eta_j-\eta_\ell)
\ee 
where $\mu\neq 0$ and $1\leq m\neq j\neq\ell\neq m\leq k$.
\vspace{1ex}

 If $\eta_1$ is a principal normal of multiplicity at least two, it 
follows from \eqref{newsame}, \eqref{hat}, \eqref{same2} and \eqref{what} that
$$
(n-2)\<\hat\eta_1,\hat\eta_j\>
+\|\hat\eta_j\|^2-(n-1)\|\hat\eta_1\|^2=0,\;\;
2\leq j\leq k,
$$
which is equivalent to
\be\label{hat5}
\|2\hat\eta_j+(n-2)\hat\eta_1\|=n\|\hat\eta_1\|.
\ee
If $\eta_i\neq\eta_j$ are principal normals of multiplicity one,
we have from \eqref{newsame} and \eqref{what} that
\be\label{hat4}
\|\hat\eta_i\|^2+(n-2)\<\hat\eta_i,\hat\eta_j\>+\|\hat\eta_j\|^2 
=n\|H\|^2-\frac{\tau}{n-1}\cdot
\ee

By Lemma \ref{unique} there is at most one  principal 
normal $\eta_1$ of multiplicity at least two. Suppose first that 
this is the case. Due to \eqref{hat5} the vectors 
$$
\b_j=2\hat\eta_j+(n-2)\hat \eta_1,\;\;2\leq j\leq k,
$$
satisfy
\be\label{samel}
\|\b_j\|=n\|\hat\eta_1\|,\;\;2\leq j\leq k.
\ee

\noindent Case $(i)$: If $j,m,\ell\geq 2$, we have from \eqref{hipo} that
$$
(1-\mu)\b_j=\b_m-\mu\b_\ell
$$
if $j\neq m\neq\ell\neq j$. Hence
$$
\|\b_j\|^2-2\mu\|\b_j\|^2+\mu^2\|\b_j\|^2=
\|\b_m\|^2-2\mu\<\b_m,\b_\ell\>+\mu^2\|\b_\ell\|^2
$$
which gives 
$$
\|\b_j\|^2=\<\b_m,\b_\ell\>.
$$
It follows from \eqref{samel} that $\eta_m=\eta_\ell$, and 
this is a contradiction. 
\vspace{1ex}

\noindent Case $(ii)$: If $j=1$ and $m\neq\ell\geq 2$, 
we have from \eqref{hipo} that
$$
\b_m-\mu\b_\ell=n(1-\mu)\hat\eta_1.
$$
Hence 
$$
\|\b_m\|^2-2\mu\<\b_m,\b_\ell\>+\mu^2\|\b_\ell\|^2
=n^2(1-\mu)^2\|\hat\eta_1\|^2.
$$
Using \eqref{samel} we obtain $\eta_m=\eta_\ell$, 
and this is a contradiction. 
\vspace{1ex}

\noindent Case $(iii)$: If $m=1$ and $j\neq\ell\geq 2$, 
we have from \eqref{hipo} that
$$
(1-\mu)\b_j+\mu\b_\ell=n\hat\eta_1.
$$
We may assume that $\mu\neq 1$ since, otherwise, we already have
a contradiction. Hence
$$
(1-\mu)^2\|\b_j\|^2+2\mu(1-\mu)\<\b_j,\b_\ell\>+\mu^2\|\b_\ell\|^2
=n^2\|\hat\eta_1\|^2.
$$
Using \eqref{samel} and $\mu\neq 1$ we have $\eta_j=\eta_\ell$,  
and this is a contradiction.
\vspace{1ex}

Next assume that all principal normals have multiplicity one. From 
\eqref{hat4} we have
$$
\|\hat\eta_i\|^2+(n-2)\<\hat\eta_i,\hat\eta_j\>+\|\hat\eta_j\|^2=nb
$$
where $i\neq j$ and 
$$
b=\|H\|^2-\frac{\tau}{n(n-1)}\cdot
$$
This is equivalent to
\be\label{hat42}
\|\beta_i^j\|^2=4nb+n(n-4)\|\hat\eta_j\|^2
\ee
where 
$$
\beta_i^j=2\hat\eta_i+(n-2)\hat\eta_j.
$$
We have from \eqref{hipo} that  
$$
(1-\mu)\b_j^i=\b_m^i-\mu\b_\ell^i
$$
if $i\neq m\neq j\neq\ell\neq i$. Hence
$$
\|\b_j^i\|^2-2\mu\|\b_j^i\|^2+\mu^2\|\b_j^i\|^2=
\|\b_m^i\|^2-2\mu\<\b_m^i,\b_\ell^i\>+\mu^2\|\b_\ell^i\|^2.
$$
We obtain using \eqref{hat42} and $i\neq j$  that
$$
\|\b_j^i\|^2=\<\b_m^i,\b_\ell^i\>.
$$
It follows that $\eta_m=\eta_\ell$, and this is 
a contradiction. \vspace{1ex}\qed

\noindent \emph{Proof of Theorem \ref{one}:}  The case $k=1$ is trivial.
In order to conclude holonomicity it is a standard fact that it suffices 
to show that the distributions 
$E_j^\perp=\oplus^k_{i=1, i\neq j}E_i$ are integrable for 
$1\leq j\leq k$  (see [11]). 
The Codazzi equation is easily seen to yield
\be\label{c1}
\<X,Y\>\nabla_Z^\perp\eta_i=\<\nabla_XY,Z\>(\eta_i-\eta_j)
\ee
and
\be\label{c2}
\<\nabla_XV,Z\>(\eta_j-\eta_\ell)=\<\nabla_VX,Z\>(\eta_j-\eta_i)
\ee
for any $X,Y\in E_i,Z\in E_j$ and $V\in E_\ell$ where
$1\leq i\neq j\neq\ell\neq i\leq k$.

It follows from \eqref{c1} that the $E_i$'s are integrable. 
Thus, it is sufficient to argue for the case $k\geq 3$.
In fact, it suffices to show that if $X\in E_i$ and $Y\in E_j$ then 
$[X,Y]\in E_\ell^\perp$ if $i\neq j\neq\ell\neq i$. 
We have from \eqref{c2} that
$$
\<\nabla_XY,Z\>(\eta_\ell-\eta_j)=\<\nabla_YX,Z\>(\eta_\ell-\eta_i)
$$
for any $Z\in E_\ell$. Then we obtain from Lemma \ref{li} that
$$
\<\nabla_X Y,Z\>=\<\nabla_Y X,Z\>=0
$$ 
which completes the proof of holonomicity. 
Then Lemma \ref{unique} completes the proof.\qed

\begin{example}\po\label{ex} 
{\em A large family of nontrivial examples of conformally flat 
$n$-dimensional submanifolds in $\R^{n+2}$ was constructed in \cite{df}.  
This construction goes as follows. Start with two smooth spherical curves 
parametrized by arc-length
$$
\gamma_i\colon I_i\subset\R\to\Sf^{m_i}_{r_i}\subset\R^{m_i+1},\ 1\leq i\leq 2,
$$ 
where $m_1+m_2=n$ and $r_1^2+r_2^2=1$. Consider the spherical surface
parametrized by the isometric immersion 
$h\colon L^2=I_1\times I_2\to\Sf^{n+1}_1\subset\R^{n+2}$ 
defined by 
$$
h(u,v)=(\gamma_1(u),\gamma_2(v)).
$$
Then, the  $n$-dimensional  submanifold of $\R^{n+2}$ parametrized on 
the unit normal bundle $UN_hL$ of $h$ in $\Sf^{n+1}_1$ by the map 
$$
\phi(w)=h(u,v)+i_*w,
$$
where $i\colon\Sf_1^{n+1}\to\R^{n+2}$ denotes the inclusion, is 
conformally flat.  A straightforward computation shows that the 
submanifold has flat normal bundle.
}\end{example}

\section{Proof of Theorem \ref{quasi}}

We first recall from \cite{dft} or \cite{dt2} the following facts which 
can easily be proved using \eqref{equiv}.

\begin{lemma}\po\label{prop1}
Let $f\colon M^n\to\R^N$ be an isometric immersion with flat normal bundle 
and principal normal vectors $\eta_1,\dots,\eta_\ell$ at $x\in M^n$. 
Denote
$$
d=\dim\spa\{\eta_i:1\leq i\leq\ell\}
$$ 
and 
$$
S_f=\spa\{\eta_i-\eta_j:1\leq i,j\leq\ell\}.
$$
\begin{itemize}
\item[(i)] $\dim S_f\leq \ell-1$ and $d-1\leq \dim S_f\leq d$.
\item[(ii)] If $\dim S_f=d-1$ then the unit vector 
$\delta\in\spa\{\eta_i,\, 1\leq i\leq\ell\}$ orthogonal to $S_f$ is 
umbilical, that is, $A_\delta=aI$.
\end{itemize}
\end{lemma}

\noindent \emph{Proof of Theorem \ref{quasi}:}
It is well-known that $M^n$, $n\geq 4$, is conformally flat if and 
only if at any $x\in M^n$ the following holds:
\be\label{prop2}
K(X_1,X_2)+K(X_3,X_4)=K(X_1,X_3)+K(X_2,X_4)
\ee
for every quadruple of orthogonal vectors $X_1,X_2,X_3,X_4\in T_xM$.

Let $\eta_1,\dots,\eta_{n-m+1}$ be the principal normals of $f$ with 
$\eta_1$  the one of multiplicity $m\geq 2$. Choosing   $X_1,X_2\in E_1$, $X_3\in E_i$ and $X_4\in E_j$, 
we obtain from 
\eqref{prop2} that
\be\label{4}
\<\xi_i,\xi_j\>=0,\;\;\mbox{for all}\;\;2\leq i\neq j\leq n-m+1,
\ee
where 
$$
\xi_i=(\eta_1-\eta_i)/\|\eta_1-\eta_i\|,\;\; 2\leq i\leq n-m+1.
$$
It follows from Lemma \ref{prop1}  that $\dim S_f=n-m$. 
In particular $p\geq n-m$.

Observe that \eqref{4} is equivalent to 
\be\label{2}
\<\xi_i,\eta_j\>=\<\xi_i,\eta_1\>\;\;\text{for any}\;\;2\leq j\neq i\leq n-m+1.
\ee
According to Lemma \ref{prop1} we have to distinguish the following cases.
\medskip 

\noindent If $\dim S_f=d-1$ we have from part $(ii)$ that 
$$
\spa\{\eta_i,\; 1\leq i\leq n-m+1\}
=\spa\{\delta,\xi_2,\dots,\xi_{n-m+1}\}.
$$
From \eqref{2} each $A_{\xi_i}$, $2\leq i\leq n-m+1$,
has an eigenvalue of multiplicity at least $n-1$. 
\medskip

\noindent If $\dim S_f=d$, then 
$$
\spa\{\eta_i,\; 1\leq i\leq n-m+1\}=\spa\{\xi_2,\dots,\xi_{n-m+1}\},
$$
and the proof follows similarly.\qed

\section{Proof of Theorem \ref{gen}}

We first define  generalized cylinders in space forms and subsequently 
generalized cones. Notice that the latter submanifolds are also generalized 
cylinders.
\vspace{1ex}

Let $g\colon L^{n-s}\to\Q_c^N$, $1\leq s\leq n-1$, be an isometric 
immersion carrying a parallel flat normal subbundle 
$\pi\colon\mathcal{L}\subset N_gL\to L^{n-s}$ of rank $s$. 
The $s$-generalized cylinder over $g$ determined by $\mathcal{L}$
is the submanifold parametrized (at the open subset of regular points) by
the map $f\colon\mathcal{L}\to\Q_c^N$ given by 
$$
f(x,v)=\exp_{g(x)}v
$$
where $\exp$ is the exponential map of $\Q_c^N$.
\vspace{1ex} 

The following result can be found in \cite{dft} or \cite{dt2}.

\begin{proposition}\po\label{cylinders} Let $f\colon M^n\to\Q_c^N$ be 
an isometric immersion with constant index of relative nullity $\nu>0$ 
such that the conullity distribution $x\in M^n\mapsto\Delta^\perp(x)$ 
is integrable. Then $f$ is locally  a $\nu$-generalized cylinder over a leaf
$g\colon L^{n-\nu}\to\Q_c^N$ of conullity.
\end{proposition}

Let $g\colon L^n\to\Q_{\tilde c}^m$ be an isometric immersion, 
and let $i\colon \Q_{\tilde c}^m\to \Q_c^N$, $\tilde c\geq c,$
be an umbilical inclusion. Since the normal bundle of 
$\tilde g=i\circ g$ splits as 
$$
N_{\tilde g} L=i_*(N_g L)\oplus N_i\Q_{\tilde c}^m,
$$
we regard $\mathcal{L}=N_i\Q_{\tilde c}^m$ as a 
subbundle of $N_{\tilde g}L$.  The $(N-m)$-generalized cone over $g$
is the submanifold parametrized (at regular points) by
the map $f\colon\mathcal{L}\to\Q_c^N$ given by 
$$
f(x,v)=\exp_{g(x)}v
$$
where $\exp$ is the exponential map of $\Q_c^N$.
\vspace{1ex}

The following result can be found in \cite{dks} or \cite{dt2}.

\begin{proposition}\po\label{prop}
Let $f\colon M^n\to \Q_c^N$ be an isometric immersion with constant 
index of relative nullity $\nu>0$. Assume that the conullity distribution 
is umbilical. Then $f$ is locally a $\nu$-generalized cone over a leaf of
conullity $g\colon L^{n-\nu}\to\Q_{\tilde c}^{N-\nu}$ contained in an 
umbilical submanifold $\Q_{\tilde c}^{N-\nu}$ of $\Q_c^N$ with 
$\tilde{c}\geq c$.
\end{proposition}

\noindent\emph{Proof of Theorem \ref{gen}:} By Theorem \ref{one} the 
conullity distribution  is integrable. Proposition~\ref{cylinders}  
asserts that $f$ is locally an open neighborhood of a 
$\nu$-generalized cylinder over a leaf 
$g=f\circ h\colon L^{n-\nu}\to\Q_c^N$ of the 
conullity distribution, where we denote by $h\colon L^{n-\nu}\to M^n$ 
the inclusion map. 

Let $u_1,\dots,u_n$ be principal coordinates for $f$ with 
corresponding coordinate vector fields $\d_1,\dots,\d_n$ such that 
$\Delta=\spa\{\d_1,\dots,\d_\nu\}$. Since the Levi-Civita connection satisfies (see [11, Proposition 1.12])
\be\label{holo}
\<\n_{\d_i} \d_j,\d_k\>=0,\;\;1\leq i\neq j\neq k\neq i\leq n,
\ee
we have  
$$
\<\n_{\d_i}\d_j,Z\>=0
$$
for $\nu+1\leq i\neq j\leq n$ and $Z\in \Delta$. Thus
$$
\alpha_h(\tilde\d_i,\tilde\d_j)=0,\;\;\nu+1\leq i\neq j\leq n,
$$
where $h_*\tilde\d_i=\d_i,\;\nu+1\leq i\leq n$. Recalling that $f$ 
is holonomic by Theorem \ref{one}, we obtain
$$
\alpha_g(\tilde\d_i,\tilde\d_j)=f_*(\alpha_h(\tilde\d_i,\tilde\d_j))
+\alpha_f(\d_i,\d_j)=0,\;\; \nu+1\leq i\neq j\leq n,
$$
which proves that $g$ is holonomic. 

We may assume $\nu=1$ since for $\nu\geq 2$ we have from \eqref{prop2} 
choosing $X_3,X_4\in\Delta$ and the Gauss equation 
that $M^n$ has constant sectional curvature $c$. 
Let $0=\eta_1,\eta_2,\ldots,\eta_k$ be the distinct principal normals. 
 From \eqref{prop2} we have
\be\label{pna}
\<\eta_i,\eta_j\>=\lambda,\;\; 2\leq i\neq j\leq k.
\ee
Moreover, if $k<n$ and $\eta_r$ is the principal normal vector field with 
multiplicity higher than one, then
\be\label{pnaa}
\|\eta_r\|^2=\lambda.
\ee
In fact, the above follows from \eqref{prop2} by choosing the 
quadruple of orthogonal vectors so that 
$X_4\in \Delta, X_1,X_2 \in E_r$ and $X_3 \in E_j$, 
or $X_1 \in E_i, X_2 \in E_j$ and $X_3 \in E_r$. 

We may assume $\lambda\neq 0$ since for $\lambda=0$ we have from 
\eqref{pna}  and the Gauss equation that $M^n$ has 
constant sectional curvature $c$. We claim that the distribution 
$\Delta^\perp$ is umbilical. Since the claim is trivial if $k=2$, 
we assume $k\geq 3$. From \eqref{c1} we have 
$$
\n^\perp_Z\eta_i=\<\n_{Y_i}Y_i,Z\>\eta_i,\;\;2\leq i\leq k,
$$
where $Y_i\in E_i$ is of unit length and $Z\in \Delta$.
Then \eqref{pna} yields
$$
Z(\lambda)=\lambda(\<\n_{Y_i}Y_i,Z\>+\<\n_{Y_j}Y_j,Z\>),\;\; 
2\leq i\neq j\leq k,
$$ 
whereas \eqref{pnaa} gives
$$
Z(\lambda)=2\lambda\<\n_{Y_r}Y_r,Z\>.
$$
Since $n\geq 4$, we obtain 
$$
\<\n_{Y_i}Y_i,Z\>=\frac{Z(\lambda)}{2\lambda},\;\; 2\leq i\leq k.
$$
On the other hand, we have from \eqref{c2} that 
$$
\<\n_{Y_i}Y_j,Z\>=0,\;\;2\leq i\neq j\leq k,
$$
and the claim follows.

Proposition \ref{prop} now gives  that $f$ coincides locally 
with the $1$-generalized cone over a leaf of conullity 
$g\colon L^{n-1}\to \Q_{\tilde c}^{N-1}$ into an umbilical 
submanifold $i\colon\Q_{\tilde c}^{N-1}\to\Q_c^N$, 
$\tilde c\geq c$. It is easy to see that
$$
\a_f(X,Y)=i_*\a_g(X,Y)
$$
for any $X,Y\in TL$. Hence $L^{n-1}$ has sectional curvature 
$\tilde c+\lambda|_L$. 

It remains to show that $\lambda$ is constant along the 
conullity leaves. From \eqref{c1} we have
$$
\n^\perp_{X_m}\eta_i=\<\n_{Y_i}Y_i,X_m\>(\eta_i-\eta_m)
$$
where $Y_i\in E_i$ is of unit length and $X_m\in E_m$ with 
$2\leq m\neq i\leq k$. If $k>3$ then \eqref{pna} gives
$$
X_m(\lambda)=0,\;\;2\leq i\neq j\neq m\neq i\leq k.
$$
If $k=3$ we have from \eqref{pna} and \eqref{pnaa} that
$$
X_i(\lambda)=2\<\n_{X_i}^\perp \eta_r,\eta_r\>=0,\;\;2\leq i\neq r\leq 3,
$$
where $X_i\in E_i$. Since $\eta_r$ is Dupin the proof is complete.\qed

\begin{example}\po{\em Take any isometric immersion with flat normal 
bundle of an open subset of $\Q_c^n$ into $\Q_c^{n+p}$ with $p<n$. 
Then it is well-known that $\nu\geq n-p$ (see Example 1 and Corollary 1 in \cite{mo}). Hence any component of the open 
dense subset where the immersion is proper is an example of a 
generalized cylinder with constant sectional curvature.
}\end{example}

\section{Proof of Theorem \ref{two}}

To prove the theorem, we first establish a one-to-one correspondence 
of local nature between globally conformally flat holonomic submanifolds 
$f\colon M^n\to\R^N$, $n\geq 4$, and flat holonomic submanifolds 
$F\colon M_0^n\to\Lo^{N+2}$  that lie inside the light-cone $\V^{N+1}$ 
of the standard flat Lorentzian space form $\Lo^{N+2}$. 
Here  $M_0^n$ denotes $M^n$ endowed with the flat metric conformal 
to the one of $M^n$. On the other hand, we have from the results in 
\cite{dt} and \cite{dt1} that any flat submanifold $M_0^n$ in $\Lo^{N+2}$ 
admits locally an abundance of Ribaucour transformations with induced 
flat metric. Then, after restricting to the transforms that preserve 
lying in the light-cone, we obtain by means of the correspondence an 
$N$-parameter family of new conformally flat holonomic submanifolds 
in $\R^N$.
\vspace{1ex}

Let $\<\,,\,\>_*$ be a metric conformal to the one of the Riemannian
manifold $(M^n,\<\,,\,\>)$ with conformal factor 
$e^{2\omega}$, that is, 
$$
\<\,,\,\>_*=e^{2\omega}\<\,,\,\>.
$$ 
The corresponding Levi-Civita connections $\n^*$ and $\n$ 
are related by
\be\label{con}
\n^*_XY=\n_XY+Y(\omega)X+X(\omega)Y-\<X,Y\>\grad\omega
\ee
where the gradient is computed with respect to the metric $\<\,,\,\>$.
 From \cite{ku} the relation between the curvature tensors 
$R_*$ and $R$ is given by
\be\label{conformal}
R_*(X,Y)Z=R(X,Y)Z-T(X,Y)Z
\ee
where
\begin{align*}
 T(X,Y)Z
=& \left(Q(Y,Z)+\<Y,Z\>\|\grad \omega\|^2\right)X  
-\left(Q(X,Z)+\<X,Z\>\|\grad \omega\|^2\right)Y \\ 
&+\<Y,Z\>Q_0X-\<X,Z\>Q_0Y,
\end{align*}
$$
Q(X,Y)=\hess\omega(X,Y)-X(\omega)Y(\omega),\;\;\;
Q_0X=\n_X \grad\omega-X(\omega)\grad\omega
$$
and everything is computed with respect to the metric $\<\,,\,\>$.
\vspace{1ex}

In the sequel $(M^n,\<\,,\,\>)$ stands for a globally  conformally 
flat manifold, that is, we have globally that 
$$
\<\,,\,\>=e^{2\omega}\<\,,\,\>_0
$$ 
where $\omega\in C^\infty(M)$ and $\<\,,\,\>_0$ is a flat metric on $M^n$.  

\begin{lemma}\po\label{Q1} Let $f\colon M^n\to\R^N$, $n\geq 3$, be a proper 
isometric immersion with flat normal bundle.  Then
\be\label{15}
Q(X,Z)=0
\ee
for $X\in E_i$ and any $Z\in TM$ such that $Z\perp X$. 
\end{lemma}

\proof Let $R^f$ denote the curvature tensor with respect to the metric 
$\<\,,\,\>$. Since  $\<\,,\,\>_0$ is flat, we have from \eqref{conformal} 
that
\be\label{1}
R^f(X,Y)Z=-T(X,Y)Z
\ee
for any $X,Y,Z\in TM$.

Assume that either $Z\in E_i$ or that $Z\in E_j$ with $j\neq i$. In either case we have $\alpha_f(X,Z)= \alpha_f(Y,Z)=0$ for $Y\perp\spa\{X,Z\}$. 
The Gauss equation gives
$$
\<R^f(X,Y)Z,W\>=\<\alpha_f(X,W),\alpha_f(Y,Z)\>
-\<\alpha_f(X,Z),\alpha_f(Y,W)\>=0
$$
for any $W\in TM$.  Using $Z\perp X$ we have from \eqref{1} that
$$
0=T(X,Y)Z=Q(Y,Z)X-Q(X,Z)Y,
$$
and \eqref{15} follows. \qed

\begin{lemma}\po\label{Q3} Let $f\colon M^n\to\R^N$ be a proper isometric 
immersion with flat normal bundle.
If $\dim E_1\geq 2$, then 
$$
Q(Z,Z)=-\frac{1}{2}\big(\|\eta_1\|^2+e^{-4\omega}\|\grad\omega\|^2\big) 
$$
for any $Z\in E_1$ with $\|Z\|=1$, where the gradient is taken with respect to the flat metric.
\end{lemma}

\proof We have from \eqref{1} that 
$$
\<R^f(X,Y)Y,X\>=-\<T(X,Y)Y,X\>
$$
where $X,Y\in E_1$ satisfy $X\perp Y$ and $\|X\|=\|Y\|=1$.
The Gauss equation yields
$$
\<R^f(X,Y)Y,X\>=\|\eta_1\|^2.
$$
Since
$$
\<T(X,Y)Y,X\>=Q(X,X)+Q(Y,Y)+e^{-4\omega}\|\grad\omega\|^2,
$$
we obtain 
\be\label{14}
Q(X,X)+Q(Y,Y)=-\|\eta_1\|^2-e^{-4\omega}\|\grad\omega\|^2.
\ee
Setting 
$$
Z=\frac{1}{\sqrt{2}}(X+Y)
$$ 
we have that 
$$
Q(Z,Z)=\frac{1}{2}(Q(X,X)+2Q(X,Y)+Q(Y,Y)),
$$
and the proof follows from \eqref{15} and \eqref{14}. 
\vspace{1ex}\qed

Let $(\Lo^{N+2},\<\!\<\,,\,\>\!\>)$ denote the standard flat Lorentzian space 
form. The \emph{light-cone} $\V^{N+1}$ of $\Lo^{N+2}$ is  
one of the connected components of the set of vectors 
$$
\{v\in\Lo^{N+2}\smallsetminus\{0\}:\<\!\<v,v\>\!\>=0\}
$$
endowed with a degenerate metric induced from $\Lo^{N+2}$.

Let $F\colon M^n\to\Lo^{N+2}$ be an isometric immersion  
of a Riemannian manifold $(M^n,\<\,,\,\>)$ that lies inside $\V^{N+1}$. 
Taking derivatives of $\<\!\<F,F\>\!\>=0$ yields that the position 
vector field $F$ is a parallel normal vector field to $F$ such that 
the second fundamental form satisfies
\be\label{umbilic}
\<\!\<\alpha_F(X,Y),F\>\!\>=-\<X,Y\>
\ee
for all $X,Y\in TM$.

Fix $w\in\V^{N+1}$.  We have from \cite{to} that the subset 
of the light-cone 
$$
\E^N=\E_w^N=\{v\in\V^{N+1}:\<\!\<v,w\>\!\>=1\}
$$
is a model for $\R^N$.  In fact, fix 
$v\in\E^N$ and let $A\colon \R^N\to(\spa\{v,w\})^\perp\subset\Lo^{N+2}$
be a linear isometry.
Then the map $\Psi=\Psi_{v,w,A}\colon\R^N\to\Lo^{N+2}$ given by 
$$
\Psi(x)=v+Ax-\frac{1}{2}\|x\|_{\R^N}^2 w
$$
is an isometric embedding such that $\Psi(\R^N)=\E^N$. Moreover, 
the normal bundle is $N_\Psi\R^N=\spa\{\Psi,w\}$ and the second 
fundamental form is given by
$$
\alpha_\Psi(X,Y)=-\<X,Y\>w
$$
for all $X,Y\in T\R^N$.
\vspace{1ex}

The map $F\colon M_0^n\to\Lo^{N+2}$ laying inside  $\V^{N+1}$ given by 
$$
F=e^{-\omega}\Psi\circ f
$$ 
is an isometric immersion of $M^n_0=(M^n,\<\,,\,\>_0)$ that was
called  in \cite{nt} the \emph{flat lift} of $f$ (see Section 9.4 in \cite{dt2}). Clearly $F$ has flat normal bundle.

\begin{lemma}\po\label{sfflorentz}
The second fundamental form of the flat lift $F$ of $f$ is given by 
$$
\alpha_F(X,Y)=-Q(X,Y)F+e^{-\omega}\Psi_*\left(\alpha_f(X,Y)
-\<X,Y\>_0f_*\grad\omega\right)-e^\omega\<X,Y\>_0 w
$$
where everything is computed with respect to the flat metric.
\end{lemma}

\proof Let  $\bar\n$, $\tilde\n$ and $\n$ denote the Levi-Civita 
connections of the metrics $\<\!\<\,,\,\>\!\>$, $\<\,,\,\>$ and  
$\<\,,\,\>_0$, respectively. Then
$$
\alpha_F(X,Y)=\bar\n_X F_*Y-F_*\n_XY.
$$
Since
$$
F_*Y=\bar\n_Y F=-e^{-\omega}Y(\omega)\Psi\circ f
+e^{-\omega}(\Psi\circ f)_*Y, 
$$
we have 
$$
\bar\n_X F_*Y=(X(\omega)Y(\omega)-XY(\omega))F
-e^{-\omega}(\Psi\circ f)_*\big(Y(\omega)X+X(\omega)Y\big)
+e^{-\omega} \bar\n_X(\Psi\circ f)_*Y.
$$
On the other hand,
\begin{align*}
 \bar\n_X(\Psi\circ f)_*Y
&=\Psi_*\n^0_Xf_*Y+\alpha_\Psi(f_*X,f_*Y)\\
&= (\Psi\circ f)_*\tilde\n_XY+\Psi_*\alpha_f(X,Y)-\<X,Y\>w
\end{align*}
where $\n^0$ is the Levi-Civita connection in $\R^N$. 
Using \eqref{con} we obtain
\begin{align*}
 \bar\n_X F_*Y 
=&\,(X(\omega)Y(\omega)-XY(\omega))F
-e^{-\omega}(\Psi\circ f)_*\big(Y(\omega)X
+X(\omega)Y-\tilde\n_{X} Y\big)\\
&+e^{-\omega}\big(\Psi_*\alpha_f(X,Y)
-\<X,Y\>w\big)\\ 
=&\,(X(\omega)Y(\omega)-XY(\omega))F
-e^{-\omega}(\Psi\circ f)_*\big(\<X,Y\>_0\grad\omega
-\n_XY\big)\\
&+e^{-\omega}\big(\Psi_*\alpha_f(X,Y)-\<X,Y\>w\big).
\end{align*}
Since
$$
F_*(\n_XY)=-\n_XY(\omega)F+e^{-\omega}(\Psi\circ f)_*\n_XY,
$$
we obtain
\begin{align*}
\alpha_F(X,Y)
=&\,\bar\n_X F_*Y-F_*\n_XY\\
=&\,(X(\omega)Y(\omega)-\hess\omega (X,Y))F
-e^{-\omega}\Psi_*\<X,Y\>_0f_*\grad\omega\\
& +e^{-\omega} \big(\Psi_*\alpha_f(X,Y)
-\<X,Y\>w\big)\\
=&\,-\;Q(X,Y)F+e^{-\omega} \Psi_*\left(\alpha_f(X,Y)
-\<X,Y\>_0f_*\grad\omega\right)-e^\omega\<X,Y\>_0w,
\end{align*}
and this concludes the proof.\qed

\begin{proposition}\po\label{lifts} 
Let $M^n$, $n\geq 4$, be a globally conformally flat Riemannian manifold
and let $f\colon M^n\to\R^N$ be a proper isometric immersion with flat 
normal bundle. Then the flat lift 
$F\colon M_0^n\to\Lo^{N+2}$ of $f$ is locally holonomic 
and proper with respect to the same principal coordinates. 

Let $F\colon M_0^n\to\V^{N+1}\subset\Lo^{N+2}$, $n\geq 4$, be a proper isometric 
immersion with flat normal bundle of a flat  Riemannian manifold. 
Then $F$ is the flat lift of an isometric immersion $f\colon M^n\to\R^N$ of 
a globally conformally flat Riemannian manifold that is locally holonomic 
and proper.
\end{proposition}

\proof By Theorem \ref{one}, there is a local coordinate system 
$(u_1,\dots,u_n)$ such that 
$$
\alpha_f(X_i,X_j)=0,\;\;1\leq i\neq j\leq n,
$$
where $X_i=\d_i/\|\d_i\|$ for $1\leq i\leq n$. From 
Lemma \ref{Q1} and Lemma \ref{sfflorentz} we obtain 
$$
\alpha_F(X_i,X_j)=0,\;\;1\leq i\neq j\leq n,
$$
hence $F$ is holonomic with respect to the same coordinate system. 

We have from Lemma \ref{sfflorentz} that
$$
\alpha_F(X_i,X_i)=-Q(X_i,X_i)F +e^{-\omega} \Psi_*(\alpha_f(X_i,X_i)
-e^{-2\omega}f_*\grad\omega)-e^{-\omega}w,\;\; 1\leq i\leq n.
$$
Since by Theorem \ref{one} at most one of the principal normals 
of $f$ has multiplicity larger than one, 
it now follows from Lemma \ref{Q3} that also $F$ is proper.
\vspace{1ex}

We now prove the second statement. The map $f\colon M^n\to \R^N$ given by
$$
\Psi\circ f=\frac{1}{\<\!\<F,w\>\!\>}F
$$
(for appropriate $w$) is an isometric immersion with respect to the 
conformally flat metric 
$$
\<\,,\,\>=\frac{1}{\<\!\<F,w\>\!\>^2}\<\,,\,\>_0.
$$
Let $\{Y_1,\dots,Y_n\}$ be an orthonormal tangent base at  $x\in M_0^n$
such that
$$
\alpha_F(Y_i,Y_j)=0,\;\; i\neq j.
$$
It follows from Lemma \ref{sfflorentz} that
$$
\Psi_*\alpha_f(Y_i,Y_j)=\frac{1}{|\<\!\<F,w\>\!\>|}Q(Y_i,Y_j)F
$$
and
$$
\alpha_F(Y_i,Y_i)=-Q(Y_i,Y_i)F+|\<\!\<F,w\>\!\>|\Psi_*\left(\alpha_f(Y_i,Y_i)
+f_*\grad\log|\<\!\<F,w\>\!\>|\right)-\frac{w}{|\<\!\<F,w\>\!\>|}
$$
for $1\leq i\neq j\leq n$. Taking norms in the first equation yields
$$
\alpha_f(Y_i,Y_j)=0,\  1\leq i\neq j \leq n.
$$
Hence $f$ has flat normal bundle. Moreover,  $f$ is proper since by the 
second equation and Lemma \ref{Q3} there is a one to one correspondence 
between the principal normals of $f$ and the ones of $F$. Finally, we 
have from Theorem \ref{one} that  $f$ is holonomic. \vspace{1ex}\qed

Let $F\colon M^n\to\Lo^{N+2}$ be an isometric immersion of a 
Riemannian manifold.  According to \cite[Theorem 17]{dt1}
any Ribaucour transformation $\tilde{F}\colon M^n\to\Lo^{N+2}$ 
of $F$ is of the form
$$
\tilde{F}=F-2\nu\varphi\F
$$
where $\F=F_*\grad\varphi+\beta$ and $\nu^{-1}=\<\!\<\F,\F\>\!\>$.
Moreover, the function  $\varphi\in C^\infty(M)$ and the vector field 
$\beta\in N_FM$ satisfy the condition
\be\label{cond}
\alpha_F(\grad\varphi,X)+\nap_X\beta=0
\ee
for any $X\in TM$. Notice that $(\varphi+c,\beta)$
also satisfies \eqref{cond} for any $c\in\R$.
\vspace{1ex}

Now assume that $F(M)\subset\V^{N+1}\subset\Lo^{N+2}$. Then 
\be\label{c}
\<\!\<F,\beta\>\!\>=\varphi + c\;\;\mbox{where}\;\; c\in\R.
\ee
In fact, using \eqref{umbilic} and \eqref{cond} it follows that
$$
X\<\!\<F,\beta\>\!\>=\<\!\<F,\nap_X\beta\>\!\>
=-\<\!\<F,\alpha_F(\grad\varphi,X)\>\!\>=X(\varphi)
$$
for any $X\in TM$.

\begin{proposition}\po Assume that $F(M)$ lies inside the light-cone 
$\V^{N+1}\subset\Lo^{N+2}$. Then the same holds for any Ribaucour 
transformation $\tilde{F}$ of $F$ for which $c=0$ in \eqref{c}.
\end{proposition}

\proof 
We have that 
$$
\<\!\<\tilde F,\tilde F\>\!\>=-4\nu\varphi\<\!\<F,\F\>\!\>
+4\nu^2\varphi^2\<\!\<\F,\F\>\!\>
=4\nu\varphi\left(\varphi-\<\!\<F,\beta\>\!\>\right)=0,
$$
and the proof follows.\vspace{1ex}\qed

Now assume that $M^n=M_0^n$ is flat and that $F\colon M_0^n\to\Lo^{N+2}$
is holonomic. It follows from \cite[Theorem 13]{dt} that the set of all
Ribaucour transformations of $F$ that preserve flatness and are holonomic 
with respect to the same principal coordinates depends on $N+1$ arbitrary
constants.
Thus, if $F(M)\subset\V^{N+1}$ it follows from the above that the family 
of Ribaucour transformations that, in addition, remain in the light-cone 
depends on $N$ parameters.
\vspace{2ex}

\noindent\emph{Proof of Theorem \ref{two}:} The proof now follows easily
using Proposition \ref{lifts}.\qed 

\begin{remark}\po {\em At least generically, the conformally flat metric 
of an element in the above family is not conformal to the original metric 
of $M^n$; see \cite[Theorem $20$]{to}.}
\end{remark}

\noindent Marcos Dajczer\\
IMPA -- Estrada Dona Castorina, 110\\
22460--320, Rio de Janeiro -- Brazil\\
e-mail: marcos@impa.br
\bigskip

\noindent Christos-Raent Onti\\
IMPA -- Estrada Dona Castorina, 110\\
22460--320, Rio de Janeiro -- Brazil\\
e-mail: christos.onti@impa.br
\bigskip

\noindent Theodoros Vlachos\\
University of Ioannina \\
Department of Mathematics\\
Ioannina -- Greece\\
e-mail: tvlachos@uoi.gr


\begin{thebibliography}{lll}

\bibitem{ca0} Cartan, E., 
\emph{La d\'eformation des hypersurfaces dans l'espace conforme r\'eel  
a $n\ge 5$ dimensions}, 
Bull. Soc. Math. France  \textbf{45} (1917), 57--121.

\bibitem{ca} Cartan, E., 
\emph{Sur les vari\'et\'es de courbure constante d'un espace 
euclidien ou non-euclidien},
Bull. Soc. Math. France \textbf{47} (1919), 125--160; \textbf{48}
(1920), 132--208.

\bibitem{bycv} Chen, B.-Y. and Verstraelen, L.,
\emph{A characterization of totally quasiumbilical submanifolds and its
applications},
Boll. Un. Mat. Ital. \textbf{14} (1977), 49--57 and 634.

\bibitem{df}  Dajczer, M. and Florit, L.,
\emph{On conformally flat submanifolds},
Comm. Anal. Geom. \textbf{4} (1996), 261--284.

\bibitem{dft}  Dajczer, M., Florit, L. and Tojeiro, R.,
\emph{Reducibility of Dupin submanifolds},
Illinois J. Math.  \textbf{49} (2005), 759--791.

\bibitem{dks} Dajczer, M.,  Kasioumis, Th., Savas-Halilaj, A. and
Vlachos, Th., \emph{Complete minimal submanifolds with nullity in the 
hyperbolic space}, J. Geom. Anal. \textbf{29} (2019), 413--427.

\bibitem{dov} Dajczer, M., Onti,  C.-R. and Vlachos, Th., 
\emph{Einstein submanifolds with flat normal bundle in space forms
are holonomic}, Proc. Amer. Math. Soc. \textbf{146} (2018), 4035--4038.

\bibitem{dt0}  Dajczer, M. and Tojeiro, R.,
\emph{Isometric immersions and the generalized Laplace and 
Elliptic Sinh-Gordon equations},  
J. Reine Angew. Math. \textbf{467} (1995), 109--147. 

\bibitem{dt}  Dajczer, M. and Tojeiro, R.,
\emph{An extension of the classical Ribaucour transformation},
Proc. London Math. Soc. \textbf{85} (2002), 211--232.

\bibitem{dt1}  Dajczer, M. and Tojeiro, R.,
\emph{Commuting Codazzi tensors and the Ribaucour transformation
for submanifolds},
Results Math. \textbf{44} (2003), 258--278.

\bibitem{dt2}  Dajczer, M. and Tojeiro, R.,
``Submanifold theory beyond an introduction".
Series: Universitext. Springer, 2019.

\bibitem{ku} Kulkarni, R., 
\emph{Curvature and metric}, 
Ann. of Math. \textbf{91} (1970), 311--331. 

\bibitem{mo0} Moore, J., 
\emph{Conformally flat submanifolds of Euclidean space}, 
Math. Ann.  \textbf{225} (1977), 89--97.

\bibitem{mo} Moore, J.,
\emph{Submanifolds of constant positive curvature I},
Duke Math. J. \textbf{44} (1977), 449--484.

\bibitem{momo} Moore, J. and  Morvan, J.-M.,
\emph{Conformally flat submanifolds of codimension four},
C. R. Acad. Sci. Paris S\'er. \textbf{287} (1978), A655--A657.

\bibitem{re} Reckziegel, H.,
\emph{Kr\"ummungsfl\"achen von isometrischen Immersionen in R\"aume 
konstante Kr\"ummung}, 
Math. Ann. \textbf{223} (1976), 169--181.

\bibitem{nt}  Donaldson, N. and Terng, C.,
\emph{Conformally flat submanifolds in spheres and integrable systems},
Tohoku Math. J. \textbf{63} (2011), 277--302.

\bibitem{to} Tojeiro, R.,
\emph{Isothermic submanifolds of Euclidean space},
J. Reine Angew. Math. \textbf{598} (2006), 1--24.
\end{thebibliography}
\end{document}